\RequirePackage[leqno]{amsmath}
\documentclass[leqno]{amsart}
\usepackage{pifont}
\usepackage{amsmath}
\usepackage{amsmath,amsthm}
\usepackage[utf8]{inputenc}
\usepackage{color}
\usepackage{amsfonts}
\usepackage{enumitem}
\setlist[itemize]{topsep=0pt,after=\vspace{1.5\baselineskip}}
 \usepackage[colorlinks=true]{hyperref}
\usepackage{empheq}
\usepackage[T1]{fontenc}

\def\R{\mathbb R}  
\newtheorem{theorem}{Theorem}[section]

\newtheorem{lemma}[theorem]{Lemma}

\newtheorem{taggedassumptionx}{Assumption}
\newenvironment{taggedassumption}[1]
 {\renewcommand\thetaggedassumptionx{#1}\taggedassumptionx}
 {\endtaggedassumptionx}

 \newtheorem{taggedefinitionx}{Definition}
\newenvironment{taggedefinition}[1]
 {\renewcommand\thetaggedefinitionx{#1}\taggedefinitionx}
 {\endtaggedefinitionx}

\newtheorem{remark}{Remark}
\def\nx{{\bf x}}

\def\io{\int_\Omega} \def\iob{\int_{\partial\Omega}} 

\newcommand{\RomanNumeralCaps}[1]
    {\MakeUppercase{\romannumeral #1}}

\title[On the lifespan of solutions to a non-local porous medium problem] 
      {On the lifespan of classical solutions to a non-local porous medium problem with nonlinear boundary conditions}
\author[M. Marras, N. Pintus and G. Viglialoro]{}
\keywords{Non-local reaction-diffusion problems, porous medium equations, global existence, blow-up.\\
\textit{$^\dagger$Corresponding author}: giuseppe.viglialoro@unica.it}

\begin{document}
\maketitle

\centerline{\scshape Monica Marras, Nicola Pintus \and Giuseppe Viglialoro$^{{\dagger}}$}
\medskip
{
 
 \centerline{Dipartimento di Matematica e Informatica}
  \centerline{Universit\`{a} di Cagliari}
  \centerline{Viale Merello 92, 09123, Cagliari (Italia)}
  
}

\bigskip
\begin{abstract}
In this paper we analyze the porous medium equation 
\begin{equation}\label{ProblemAbstract}
\tag{$\Diamond$}
u_t=\Delta u^m + a\io u^p-b u^q -c\lvert\nabla\sqrt{u}\rvert^2 \quad \textrm{in}\quad  \Omega \times I,
\end{equation}
where $\Omega$ is a bounded and smooth domain of $\R^N$, with $N\geq 1$,  and $I= [0,t^*)$ is the maximal interval of existence for $u$. The constants $a,b,c$ are positive, $m,p,q$ proper real numbers larger than 1 and  the equation is complemented with nonlinear boundary conditions involving the outward normal derivative of $u$. Under some hypothesis on the data, including intrinsic relations between $m,p$ and $q$, and assuming that for some positive and sufficiently regular function $u_0(\nx)$ the Initial Boundary Value Problem (IBVP) associated to \eqref{ProblemAbstract} possesses a positive  classical solution $u=u(\nx,t)$ on $\Omega \times I$:
\begin{itemize}
\item [$\triangleright$] when $p>q$ and in 2- and 3-dimensional domains, we determine a \textit{lower bound of} $t^*$ for those $u$ becoming unbounded in $L^{m(p-1)}(\Omega)$ at such $t^*$;
\item [$\triangleright$] when $p<q$ and in $N$-dimensional settings,  we establish a \textit{global existence criterion} for $u$.
\end{itemize}
\end{abstract}
\section{Introduction, state of the art and motivations}\label{IntroductionSection} 
Reaction-diffusion equations are commonly employed to model several natural phenomena appearing in various physical, chemical and biological applications. This paper is devoted to this specific reaction-diffusion problem 
\begin{equation}\label{General_Problem}
\begin{cases}
u_t=  \Delta u^m+ a\int_{\Omega}u^p-bu^q-c\lvert\nabla\sqrt{u}\rvert^2 &  \textrm{in} \quad \Omega \times I, \\
\nabla u \cdot\boldsymbol\nu=u_{\boldsymbol\nu}=g\left(u\right) & \textrm{on} \quad \partial \Omega \times I, \\ 
u({\bf x},0)=u_0({\bf x}) & {\bf x} \in\Omega, 
\end{cases}
\end{equation}
where $u=u({\bf x},t)$ is defined in the cylinder $\Omega\times I$, $\Omega$ being a bounded and smooth domain of $\R^N$ ($N\geq 1$) with regular boundary $\partial \Omega$, and $I= [0,t^*)$ the maximal interval of existence for the solution $u$. Further, ${\bf \boldsymbol\nu}=(\nu_1,\ldots,\nu_N)$ stands for the outward normal unit vector to the boundary $\partial \Omega$ and  $u_{ \boldsymbol\nu}$ is the normal derivative of $u$. Additionally, $u_0:=u_0({\bf x})>0$ and $g=g(\xi)\geq 0$ are appropriate and sufficiently smooth functions of, respectively, ${\bf x} \in \Omega$ and $\xi\geq 0$, $a,b,c$ arbitrarily positive constants, whilst $m,p,q>1$ are proper reals related one to the other.

Beyond problems arising in the mathematical models for gas or fluid flow in porous media (see \cite{Aronson1986} and \cite{vazquez2007porous}), the formulation in \eqref{General_Problem} also describes (see \cite{GURTIN197735} and \cite{Souplet_Gradient}) the evolution of some biological population $u$ (cells, bacteria, etc.) at some time $t$ and position ${\bf x}$ which lives in a certain domain $\Omega$ and whose growth is induced by the law $a\int_{\Omega}u^p-bu^q-c\lvert\nabla\sqrt{u}\rvert^2$ (positive addends essentially represent sources which increase the energy of the system and stimulate the occurrence of an uncontrolled increasing of $u$ through the time, whereas negative ones have a damping/absorption effect, absorb the energy and, so, contrast the power of the source terms); precisely, the term $\Delta u^m$ idealizes the spread of the population, the parameter $m>1$ indicating the speed of propagation/displacement, the non-local term $+a\int_{\Omega}u^p$ the births of the species  and $-bu^q-c\lvert\nabla\sqrt{u}\rvert^2$ counts, respectively, its natural and the accidental deaths. Moreover, the assumption $u_{\boldsymbol\nu}=g\left(u\right)\geq 0$ on the boundary virtually models an incoming flux of the population $u$; in a particular way, since realistically to a low (high) concentration of $u$ on the boundary corresponds a low (high) incoming flux,  $g$ will be taken from a set including increasing functions. (For interested readers, we mention the contributions \cite{AiTsEdYaMi, MarrasVernierVigliaWithm,MarrasVigliaMatNac, ViglialoroDifferentialIntegralEquations,ViglialoroBoundnessVeryWeak,
ViglialoroWolleyDCDS,ViglialoroWolleyMMAS2017}, where reactions terms similar to that in \eqref{General_Problem} have been also employed in chemotaxis models.)

Coming back to the mathematical analysis of our problem, some general results concerning existence of local or global solutions (i.e., respectively, $I=[0,t^*)$, $t^*$ finite or $I=[0,\infty)$), have been already studied in the literature for a class of reaction-diffusion models, where the first equation of \eqref{General_Problem} reads $\tau u_t=\nabla \cdot A({\bf x},t,u,\nabla u)+B({\bf x},t,u,\nabla u)$, $A$ and $B$ verifying some standard ellipticity behaviors as well as growth assumptions (we refer, for instance, to \cite{Kielhofer1974,Krylov,LSUBookInequality,QuittnerSouplet2007superlinear} for the case $\tau=1$ and  \cite{PucciSerrin} for $\tau=0$). In particular, as to the questions concerning existence of classical solutions to the previous equation and/or their nonnegativity (through applications of maximum principles), some results for the case $\tau=1$ can preserved if  so called ``non-degenerate'' data $u_0$ are considered (see \cite[$\S$3.1 and $\S$3.2]{vazquez2007porous}).

Despite a deep research, to the best of our knowledge the problems on the existence and regularity of solutions to \eqref{General_Problem} are not directly indicated in the present literature. For this reason, in this investigation we abstain from such an analysis, but rather we follow the same approach used in largely cited papers (see, for instance, \cite{PaynPhilSchaefer,PAYNE_Shaefer_2007,PaynePhiProytc,ShaeferWithoutGradient,ShaeferExistClassicaPorous} and references therein) where \textit{nonnegative} \textit{classical solutions} are a priori assumed to exist for a period of time but, also, the may become unbounded at some finite time $t^*$. In particular, in \cite[$\S$1]{vazquez2007porous} a discussion on the \textit{Porous Medium Equation}, $u_t=\Delta u^m$, and the \textit{Signed Porous Medium Equation}, $u_t=\Delta (|u|^{m-1}u)$, is carried out: in agreement with our purpose, we indicate that the default setting of the first case includes only nonnegative solutions.

As to well established results, there exists an important number of papers concerning variants of the IBVP \eqref{General_Problem}, some of which dealing with properties of classical solutions: global and/or local existence, lower and upper bound of blow-up time, blow-up rates and/or asymptotic behavior. In particular, we collect the following results: 
\textit{{\bf \RomanNumeralCaps{1})} $m=1$, $b=c=0$ and $\int_\Omega u^p$ replaced with $u^p$, for $p>1$}.  When $\Omega$ is a bounded and smooth domain of $\R^3$ and Dirichlet boundary conditions are assigned (i.e. $u({\bf x})=0$ for ${\bf x}\in \partial \Omega$), in \cite{PAYNE_Shaefer_2007} a lower bound for the blow-up time of solutions, if blow-up occurs, is derived, and \cite{PayneSchaefer_Robin} essentially deals with blow-up and global existence questions for the same problem in the $N$-dimensional setting, with $N\geq 2$, and endowed with Robin boundary conditions (i.e. $u_{\boldsymbol\nu}=-hu$, $h>0$, on $\partial \Omega$). 
\textit{{\bf \RomanNumeralCaps{2})} $m>1$, $b=c=0$ and $\int_\Omega u^p$ replaced with $u^p$, for $p>1$}. For $\Omega=\R^N$,  $N\geq 1$, in  \cite{galaktionov_1994}, \cite{GalaktionovKurdyMikha} and  \cite{LevineTheRoleOf} it is shown that for $1<p\le m+(2/N)$ the problem has no global positive solution, whilst for $p>m+(2/N)$ there exist initial data $u_0$ emanating global solutions. When $\Omega$ is a bounded and smooth domain of $\R^N$, $N\geq 1$, and under Dirichlet boundary conditions, in \cite{Galaktionov1981_Russian} is proved that for $1<p<m$ the problem admits global solutions for all $u_0$ such that $u_0^{m-1}\in H_0^1(\Omega),$ while for $m<p<m(1+(2/N))+(2/N)$ specific initial data produce unbounded solutions.
\textit{{\bf \RomanNumeralCaps{3})} $m>1$, $a=0$ and $-b\io u^q-c\lvert \sqrt{u}\rvert^2$ replaced with $+u^p-u^\mu |\nabla u^\alpha|^q$, with $p,q,\alpha\geq  1$ and $\mu\geq 0$}. With $\Omega$ bounded and smooth in $\R^N$, $N\geq 1$, and under Dirichlet boundary conditions, in \cite{AndreuEtAl} the authors treat the existence of the so called \textit{admissible solutions} and  show that they are globally bounded  if   $p<\mu + mq$ or $m< p=  \mu + mq$, as well as  the existence of blowing up admissible solutions, under the complementary
condition $1\leq \mu  + mq < p$. 
\textit{{\bf \RomanNumeralCaps{4})}} $m>1$, $b=c=0$ and $+a \int_\Omega u^p$ replaced with $+u^p \int_\Omega u^q$, with $p,q>1$. In a bounded and smooth domain $\Omega$ of $\mathbb{R}^3$ and under Robin boundary conditions, a lower bound for the blow-up time  if the solution blows up is determined under the assumption $p+q>m>1$ whilst conditions which ensure that the blow-up does not occur are also presented if $p+q\leq m$ (see \cite{LIURobinNonlocal}). 
\textit{{\bf \RomanNumeralCaps{5})}} $m=1$ and $c=0$, with $p,q>1$. In a bounded and smooth domain $\Omega$ of $\mathbb{R}^3$ and under proper nonlinear boundary conditions, in \cite{LIUNonLocalNoNlinearBC} a lower bound for the blow-up time  if the solution blows up is determined under the assumption $p>q.$  
\textit{{\bf \RomanNumeralCaps{6})}} $m>1$, $b=c=0$ $and$ $p\geq 0$.  In a bounded and smooth domain $\Omega$ of $\mathbb{R}^N$,  $N\geq 1$, and under Dirichlet boundary conditions, classical nonnegative solutions which are global or blow up in finite time are derived for any $0\leq u_0\in C^{2+\alpha}(\Omega) \cap C^0(\bar\Omega)$, for some $0<\alpha<1$, provided some compatibility conditions and assumptions on the data are given (see  \cite{LiXie_ExistencePorousNonLocal}). 
\textit{{\bf \RomanNumeralCaps{7})}} $m=1$, $c=0$ $and$ $p,q\geq 1$.  In a bounded and smooth domain $\Omega$ of $\mathbb{R}^N$,  $N\geq 1$, and under various boundary conditions, in \cite{WangWangNonlocalHeat} for any compatible $0\leq u_0\in C^1(\bar \Omega) $ classical nonnegative solutions which are global are attained if $p<q$, whereas for $p>q$ blowing up ones are  detected (see also \cite{SONGNonLocal}).

Motivated by the discussion so far presented, aim of the present research is expanding the underpinning theory of the mathematical analysis of problem \eqref{General_Problem}, which is not included in the above cases. In particular, the aforementioned state of the art inspires our work, and even if we will use some ideas employed in those items to address our statements, some further derivations will be necessarily required; moreover we do not restrict to prove our main theorems but we complement the general presentation of the manuscript by means of remarks and discussions.

To be precise, our contribution includes an analysis for the maximal interval $I=[0,t^*)$ of existence for classical solutions $u$ (in the sense of the Definition \ref{DefiSolution} given in $\S$\ref{SectionPreparatory} below) to system  \eqref{General_Problem}, where $t^*$ plays the role of the unknown and obeys the following \textit{extensibility criterion} (\cite{BandleBrunnerSurvey,Kielhofer1974}):  
either $t^*=+\infty$, so that $u$ remains bounded for all $\nx\in\Omega$ and all time $t>0$ and $I=(0,\infty)$, or $t^*$ is finite (blow-up time), so that $u$ exists only in $I=(0,t^*)$ and $\lVert u(\cdot,t)\rVert_{L^\infty(\Omega)} \nearrow +\infty$ as $t \searrow t^{*}$.

Thereof, we prove three theorems which provide its estimates or its precise value; they are discussed in details in $\S$\ref{mainTheoremsSection}, whilst now they are briefly summarized as follows:
\begin{itemize}
\item [$\bullet$] \textit{Lower bound of $t^*$ in $\R^3$ and $\R^2$: Theorem \ref{TheoremBlowUpDifferencePower} and Theorem \ref{TheoremBlowUpDifferencePower2}}. If for $p>q>\frac{3}{2}$, $2<m<\frac{8}{5-p}$ with $\frac{3}{2}<p<5$, $m>2$ with $p\geq 5$, $g(u)$ behaving as $u^{\frac{m(p-1)}{4}-m+2}$ and $u_0$ sufficiently regular problem \eqref{General_Problem} admits a positive classical solution $u$ which becomes unbounded in $L^{m(p-1)}(\Omega)$ at some finite time $t^{*}$, then there exists a computable $T>0$ such that $t^{*}\geq T$.
\item [$\bullet$] \textit{Criterion for global existence in $\R^N$, $N\geq 1$: Theorem \ref{TheoremGlobalferencePower}.} If for $q>p>m>1$, $2p<m+q$, $g(u)$ behaving as $u^{p-m+1}$ and $u_0$ sufficiently regular problem \eqref{General_Problem} admits a positive classical solution $u$, then holds that $t^*=\infty$.
\end{itemize}
\section{Assumptions, definitions and preparatory lemmas}\label{SectionPreparatory}
In this section we fix crucial hypothesis and lemmas which will be considered through the paper in the proofs of the main theorems. This preparatory material is herein presented according to our purposes.
\begin{taggedassumption}{$\mathcal{A}$}\label{Assumption1}
For any $N\geq 1$, $\Omega$ is a bounded and smooth domain of $\mathbb{R}^N$, star-shaped, convex in two orthogonal directions and such that, for some origin inside ${\bf x}_0$, its geometry is characterized by 
\begin{equation*}
\rho_0 := \min_{\partial \Omega} (({\bf x}-{\bf x}_0)  \cdot {\bf \boldsymbol\nu} )\;\; \textrm{and}\;\;\; d:= \max_{\overline{\Omega}} | {\bf x}-{\bf x}_0 |.
\end{equation*}
\end{taggedassumption}

\begin{taggedefinition}{$\mathcal{CS}$}\label{DefiSolution}
A \textit{classical solution} to problem \eqref{General_Problem} is a positive function $u\in C^{2,1}(\Omega \times (0,t^*))\cap C^{1,0}(\bar{\Omega} \times [0,t^*))$ which satisfies  \eqref{General_Problem}, for some $0<t^*\leq +\infty$. 
\end{taggedefinition}

\begin{taggedefinition}{$\mathcal{D}$}\label{Assumption2}
For any $p>\frac{3}{2}$ and $m>2$, let $t^*>0$ finite. We say that a nonnegative function $V\in C^0(\bar\Omega \times (0,t^*))$ blows up in $L^{m(p-1)}(\Omega)$-norm at finite time $t^*$ if 
\[\displaystyle \lim_{t\rightarrow t^*} \int_{\Omega}V^{m(p-1)}=\infty.\] 
\end{taggedefinition}
\begin{lemma}\label{lemma1} 
Let $\Omega$ be a domain satisfying Assumption \ref{Assumption1}. 
For any positive function $V\in C^1(\bar \Omega)$ and $\lambda\geq 1$, we have
\begin{equation}\label{SobolevTypeInequBoundary}
\int_{\partial \Omega} V^\lambda  \leq \frac{N}{\rho_0}\int_\Omega V^\lambda +\dfrac{d\lambda}{\rho_0}\int_\Omega V^{\lambda-1}|\nabla V|.
 \end{equation}
Moreover for every arbitrary $\epsilon>0$ we also have that:
\begin{itemize}
\item [$\triangleright$] If $N=2$
\begin{equation}\label{Inequ_v^3Dim2}
\begin{split}
\int_\Omega V^{\frac{3}{2}\lambda}   \leq & \dfrac{\sqrt{2}}{2\rho_0}\Big[\Big(\int_\Omega V^\lambda  \Big)^\frac{3}{2}+\frac{d+\rho_0}{2\epsilon^2}\Big(\int_\Omega V^\lambda \Big)^2+\frac{\left(d+\rho_0\right)\epsilon^2}{2} \int_\Omega |\nabla V^\frac{\lambda}{2}|^2 \Big]; 
\end{split}
 \end{equation}
\item [$\triangleright$] If $N=3$ 
\begin{equation}\label{Inequ_v^3}
\begin{split}
\int_\Omega V^{\frac{3}{2}\lambda}  \leq & \sqrt 2\Big[\left(\dfrac{3}{2\rho_0}\right)^\frac{3}{2}\Big(\int_\Omega V^\lambda  \Big)^\frac{3}{2}+ \frac { 1}{4 \epsilon^3}\left(1+\frac{d}{\rho_0}\right)^\frac{3}{2}\Big(\int_\Omega V^\lambda  \Big)^3\\&+\frac{3} {4}\left(1+\frac{d}{\rho_0}\right)^\frac{3}{2}\epsilon  \int_\Omega |\nabla V^\frac{\lambda}{2}|^2  \Big].
\end{split}
 \end{equation}
\end{itemize}
 \begin{proof}
 For the proofs see \cite{PPVPII} and \cite[Lemma 3.2]{ViglialoroDifferentialIntegralEquations}.
 \end{proof}
\end{lemma}
\section{Analysis and proofs of the main results}\label{mainTheoremsSection}
After the preparations in $\S$\ref{SectionPreparatory}, we are in the position to demonstrate the theorems whose general overviews were summarized in $\S$\ref{IntroductionSection}. 
\subsection{Lower bounds of the blow-up time}
The first theorem is concerned with lower bounds of the blow-up time $t^*$, through which is identified the maximal interval $I$ where solutions to system \eqref{General_Problem} are defined. We are not aware of general results indicating assumptions on the data which straightforwardly infer the existence of unbounded solutions to such a system; nevertheless, in the spirit of the results discussed in $\S\ref{IntroductionSection}$, for which blow-up may manifest for large initial data $u_0({\bf x})$ and high effects of source (coefficient $p$), or low absorption or/and diffusion (coefficient $q$ or/and $m$), we understand that also in view of the incoming flow of the population $u$, it seems reasonable to expect scenarios where unbounded solutions may appear. 
\begin{theorem}\label{TheoremBlowUpDifferencePower} 
Let $\Omega$ be a domain of $\mathbb{R}^3$ satisfying Assumption \ref{Assumption1}, $a,b,c,k>0$ and $q>\frac{3}{2}$. Moreover for any $p>q$ and 
\begin{equation}\label{RelationM-P_ValidityBlow}
\begin{cases}
m\in (2,\frac{8}{5-p})& \textrm{for} \quad \frac{3}{2}<p<5,\\
m \in (2,\infty) & \textrm{for} \quad p\geq 5,
\end{cases}
\end{equation}
let $g=g(\xi)$ be a continuous function such that $0\leq g(\xi)\leq k\xi^\beta$, for $\xi\geq 0$ and $\beta=\frac{m(p-1)}{4}-m+2>0$. If  $u$ is a classical solution, in the sense of Definition \ref{DefiSolution}, to \eqref{General_Problem} emanating from a positive initial data $u_0:=u_0({\bf x})\in C^{2+\alpha}(\Omega)\cap C^1(\bar{\Omega})$, for some $0<\alpha<1$ and such that $\nabla u_0 \cdot\boldsymbol\nu=g(u_0)$ on $\partial \Omega$, which additionally complies with Definition \ref{Assumption2}, then there exist  computable constants $c_1>0$ and $c_5>0$ such that for $\varphi(0)=\io u_0^{m(p-1)}>0$
 \begin{equation*} 
t^{*}    \geq \frac{1}{2c_1}\log \bigg(1+\frac{c_1}{c_5}\varphi^{-2}(0)\bigg)=:T.
 \end{equation*}
In particular $I=[0,t^*)\supseteq [0,T).$
\begin{proof}
If $u$ is a positive  classical solution of \eqref{General_Problem} defined in $\Omega \times (0,t^*)$ and satisfying $u_{ \boldsymbol\nu}=g(u)$ on $\partial \Omega$, by setting $s=p-1$ and using the integration by parts formula, the evolution in time of $t \mapsto \io u^{ms}$ fulfills for all $t\in (0,t^*)$
\begin{equation*}
\begin{split}
\frac{d}{d t} \io u^{ms} &= ms\io u^{ms-1}\left[\Delta u^m+a\io u^p -bu^q -c\lvert\nabla\sqrt{u}\rvert^2\right]\\&
=m^2s\iob u^{ms+m-2} u_{ \boldsymbol\nu}-m^2s(ms-1)\io u^{ms+m-3}\lvert\nabla u\rvert^2 
\\&\quad +ams\io u^{ms-1}\io u^p-bms\io u^{ms+q-1}-cms\io u^{ms-1}\lvert\nabla\sqrt{u}\rvert^2.
\end{split}
\end{equation*}
Since by the H\"{o}lder inequality  (recall $p>\frac{3}{2}$ and $m>2$) we have that 
\begin{equation*}
\io u^{ms-1}\io u^p \leq \left(\io u^{s(m+1)}\right)^{\frac{ms-1}{s(m+1)}}\left(\io u^{s(m+1)}\right)^{\frac{s+1}{s(m+1)}}\lvert\Omega\rvert\quad \textrm{for all} \quad t\in(0,t^*),
\end{equation*}
from assumption $0\leq g(u)\leq ku^{\beta}$, $u> 0$, and the pointwise identity 
$
\lvert\nabla\sqrt{u}\rvert^2=\dfrac{1}{4u}\lvert\nabla u\rvert^2\,,
$
we get neglecting the nonpositive term $-m^2s(ms-1)\io u^{ms+m-3}\lvert\nabla u\rvert^2$
\begin{equation*}
\begin{split}
\frac{d}{d t} \io u^{ms}\leq &
m^2sk\iob u^{ms+m-2+\beta}+ams\lvert\Omega\rvert\io u^{s(m+1)}\\&
-bms\io u^{ms+q-1}-\dfrac{cms}{4}\io u^{ms-2}\lvert\nabla{u}\rvert^2\quad \textrm{for all $t$}\in (0,t^*)\,.
\end{split}
\end{equation*}
For the value of $\beta$ as in our assumptions, the above relations reads
\begin{equation}\label{evolutionPhi3}
\begin{split}
\frac{d}{d t} \io u^{ms}\leq & m^2sk\iob u^{\frac{5}{4}ms}+ams\lvert\Omega\rvert\io u^{s(m+1)}\\&
-bms\io u^{ms+q-1}-\dfrac{cms}{4}\io u^{ms-2}\lvert\nabla{u}\rvert^2\quad \textrm{for all $t$}\in (0,t^*),
\end{split}
\end{equation}
so that \eqref{SobolevTypeInequBoundary} with $\lambda=\frac{5}{4}ms>1$ and $V=u$ provides
\begin{equation*}
\iob u^{\frac{5}{4}ms} \leq \dfrac{3}{\rho_0}\io u^{\frac{5}{4}ms}+\dfrac{5msd}{4\rho_0}\io u^{\frac{5}{4}ms-1}\lvert\nabla u\rvert \quad \textrm{for all $t$}\in (0,t^*),
\end{equation*}
and by plugging this gained estimate into \eqref{evolutionPhi3}, one achieves for all $t\in (0,t^*)$
\begin{equation}\label{evolutionPhi4}
\begin{split}
\frac{d}{d t} \io u^{ms}\leq & m^2sk\left[\dfrac{3}{\rho_0}\io u^{\frac{5}{4}ms}+\dfrac{5msd}{4\rho_0}\io u^{\frac{5}{4}ms-1}\lvert\nabla u\rvert\right]\\&+ams\lvert\Omega\rvert\io u^{s(m+1)}
 -bms\io u^{ms+q-1}-\dfrac{cms}{4}\io u^{ms-2}\lvert\nabla u\rvert^2.
\end{split}
\end{equation}
On the other hand, applications of the H\"{o}lder and the Young inequalities allow us to control some terms in (\ref{evolutionPhi4}). Precisely for all $t\in(0,t^*)$ holds that 
\begin{equation}\label{Inequ_v^3_2_Precedente}
\begin{split}
&
\io u^{\frac{5}{4}ms}\leq \dfrac{1}{2}\io u^{\frac{3}{2}ms} +\dfrac{1}{2}\io u^{ms},\\ &
\io u^{ms+q-1}\geq \lvert\Omega\rvert^{\frac{1-q}{ms}}\Big(\io u^{ms}\Big)^{\frac{ms+q-1}{ms}},\\&
\io u^{\frac{5}{4}ms-1}\lvert\nabla u\rvert = \dfrac{2}{ms}\io u^{\frac{3}{4}ms}\lvert\nabla u^{\frac{ms}{2}}\rvert\leq \dfrac{1}{2\epsilon_1}\io u^{\frac{3}{2}ms} +\dfrac{2\epsilon_1}{m^2s^2}\io \lvert\nabla u^{\frac{ms}{2}}\rvert^2,\\&
\io u^{s(m+1)}=\io u^{3s}u^{ms-2s}\leq \dfrac{2}{m}\io u^{\frac{3ms}{2}}+\dfrac{m-2}{m}\io u^{ms},
\end{split}
\end{equation}
which infer through \eqref{evolutionPhi4} and on the entire interval $(0,t^*)$
\begin{equation}\label{evolutionPhi5}
\begin{split}
\frac{d}{d t} \io u^{ms}\leq & m^2sk\dfrac{3}{2\rho_0}\left[\io u^{\frac{3}{2}ms}+\io u^{ms}\right]\\&
+m^2sk\dfrac{5mds}{4\rho_0}\left[\dfrac{1}{2\epsilon_1}\io u^{\frac{3}{2}ms}+\dfrac{2\epsilon_1}{m^2s^2}\io\lvert\nabla u^{\frac{ms}{2}}\rvert^2\right]\\&
+ams\lvert\Omega\rvert\left[\dfrac{2}{m}\io u^{\frac{3}{2}ms}+\dfrac{m-2}{m}\io u^{ms}\right]\\&
 -bms\lvert\Omega\rvert^{\frac{1-q}{ms}}\Big(\io u^{ms}\Big)^{\frac{ms+q-1}{ms}}
-\dfrac{cms}{4}\io u^{ms-2}\lvert\nabla u\rvert^2. 
\end{split}
\end{equation}
As to $\int_\Omega u^{\frac{3}{2}ms}$, we invoke \eqref{Inequ_v^3} to get on $(0,t^*)$
\begin{equation*}
\io u^{\frac{3}{2}ms}\leq \dfrac{3^{\frac{3}{2}}}{2\rho_0^{\frac{3}{2}}}\Big(\io u^{ms}\Big)^{\frac{3}{2}}+\sqrt{2}\left(\dfrac{d}{\rho_0}+1\right)^{\frac{3}{2}}\left[\dfrac{1}{4\epsilon_2^3}\Big(\io u^{ms}\Big)^3+\dfrac{3\epsilon_2}{4}\io \lvert\nabla u^{\frac{ms}{2}}\rvert^2\right],
\end{equation*}
so that, using the identity $u^{ms-2}\lvert\nabla u\rvert^2=\dfrac{4}{m^2s^2}\lvert\nabla u^{\frac{ms}{2}}\rvert^2$ and introducing the constants
\begin{equation*}
\begin{cases}
c_1 =a\left(m-2\right)s\lvert\Omega\rvert+\dfrac{3m^2sk}{2\rho_0},\quad 
c_2=\left(as\lvert\Omega\rvert+\dfrac{3m^2sk}{4\rho_0}+\dfrac{5m^3s^2kd}{16\rho_0\epsilon_1}\right)\dfrac{3^{\frac{3}{2}}}{\rho_0^{\frac{3}{2}}},\\
c_3=\Big(\dfrac{d}{\rho_0}+1\Big)^\frac{3}{2}\dfrac{1}{4\epsilon_2^3}\Big(2 \sqrt{2}a s|\Omega|+\dfrac{3\sqrt{2}m^2 s k}{2\rho_0}+\dfrac{5\sqrt{2}m^3s^2kd}{8\rho_0\epsilon_1}\Big),\\
c_4=\Big(\dfrac{d}{\rho_0}+1\Big)^\frac{3}{2}{\epsilon_2}\Big(\dfrac{3}{2}\sqrt{2}a s|\Omega|+\dfrac{9\sqrt{2}m^2 s k}{8\rho}+\dfrac{15\sqrt{2}m^3s^2kd}{32 \epsilon_1\rho_0}\Big)
+\dfrac{5mdk\epsilon_1}{2\rho_0}-\dfrac{c}{ms}\,,
\end{cases}
\end{equation*}
after some tedious computations, inequality \eqref{evolutionPhi5} is simplified to 
\begin{equation*}
\varphi'\leq c_1\varphi+c_2\varphi^{\frac{3}{2}}+c_3\varphi^3+c_4\io \lvert\nabla u^{\frac{ms}{2}}\rvert^2 -bms\lvert\Omega\rvert^{\frac{1-q}{ms}}\varphi^{\frac{ms+q-1}{ms}} \quad \textrm{on}\quad (0,t^*), 
\end{equation*}
where for convenience we have set $\varphi(t)=\varphi:=\io u^{ms}$ on $(0,t^*)$.

Additionally, for any fixed $0<\epsilon_1<\frac{2\rho_0 c}{5m^2sdk}$ there exists $\epsilon_2>0$ such that   $c_4\leq 0$, leading to
\begin{equation}\label{PhiPrimeCasiFinal}
\varphi'\leq c_1\varphi+c_2\varphi^{\frac{3}{2}}+c_3\varphi^3-bms\lvert\Omega\rvert^{\frac{1-q}{ms}}\varphi^{\frac{ms+q-1}{ms}}\quad \textrm{for all}\quad t\in (0,t^*).
\end{equation}
In order to obtain an \textit{explicit} estimate (see Remark \ref{RemarkExplicitTStar} below) of lower bounds for $t^*$, we do not neglect the negative addendum $-bms\lvert\Omega\rvert^{\frac{1-q}{ms}}\varphi^{\frac{ms+q-1}{ms}}$ but rather we treat it in terms of $\varphi^{\frac{3}{2}}$ and $\varphi^3$. In this sense, by using Young's inequality, we can write (recall $m>2$ and $p>q$) for any $\epsilon_3>0$
\begin{equation*}
\begin{split}
\varphi^{\frac{3}{2}}&= \big(\varphi^\frac{ms+q-1}{ms}\epsilon_3^\frac{4ms-2q+2}{3ms}\big)^\frac{3ms}{4ms-2q+2}\big(\varphi^3\epsilon_3^{-\frac{4ms-2q+2}{ms-2q+2}}\big)^\frac{ms-2q+2}{4ms-2q+2}\\ &\leq \dfrac{3ms}{4ms-2q+2}\epsilon_3^{\frac{4ms-2q+2}{3ms}}\varphi^{\frac{ms+q-1}{ms}}+\dfrac{ms-2q+2}{4ms-2q+2}\epsilon_3^{-\frac{4ms-2q+2}{ms-2q+2}}\varphi^3\quad \textrm{on}\quad  (0,t^*),
\end{split}
\end{equation*}
from which \eqref{PhiPrimeCasiFinal} is transformed into
\begin{equation*}
\begin{split}
\varphi'\leq & c_1\varphi+\Big(c_2\dfrac{3ms}{4ms-2q+2}\epsilon_3^{\frac{4ms-2q+2}{3ms}}-bms\lvert\Omega\rvert^{\frac{1-q}{ms}}\Big)\varphi^{\frac{ms+q-1}{ms}}\\&
+\left(c_2\dfrac{ms-2q+2}{4ms-2q+2}\epsilon_3^{-\frac{4ms-2q+2}{ms-2q+2}}+c_3\right)\varphi^3\quad \textrm{on}\quad  (0,t^*).
\end{split}
\end{equation*}
Finally, choosing $\epsilon_3=\left[\frac{b}{3c_2}\left(4ms-2q+2\right)\lvert\Omega\rvert^{\frac{1-q}{ms}}\right]^{\frac{3ms}{4ms-2q+2}}>0$ and, in turn, setting $c_5=c_2\frac{ms-2q+2}{4ms-2q+2}\epsilon_3^{-\frac{4ms-2q+2}{ms-2q+2}}+c_3$ we have
\begin{equation}\label{eqdiff}
\varphi(t)'\leq c_1\varphi(t)+c_5\varphi^3(t)     \quad \textrm{for all}\quad  (0,t^*).
\end{equation}
Now, since we are assuming that $\varphi(t)\nearrow \infty$ as $t\searrow t^{*}$, $\varphi(t)$ can be non decreasing, so that $\varphi(t) \geq \varphi(0)>0$ with $t\in (0,t^{*})$, or non increasing (possibly presenting oscillations), so that there exists a time $t_1$ where $\varphi(t_1)=\varphi(0)$. In any case,  $\varphi(t) \geq \varphi(0)$ for all  $t \in [t_1, t^{*})$, where $0\leq t_1<t^{*}$. Henceforth, by integrating \eqref{eqdiff} between $t_1$ and $t^{*}$, we arrive at (recall $\varphi(t_1)=\varphi(0)$) this explicit estimate for $t^*$:
\begin{equation*}\label{Lower_Dirichlet}
\begin{split}
t^*\geq \int_{\varphi(0)}^{+\infty}  \dfrac{d\tau}{c_1\tau+c_5\tau^3}
= \dfrac{1}{2 c_1}\log\left(1+\dfrac{c_1}{c_5}\varphi^{-2}(0)\right).
\end{split}
\end{equation*}
\end{proof}
\end{theorem}
\begin{remark}\label{RemarkExplicitTStar}
We point out that if in \eqref{PhiPrimeCasiFinal} we ignored the negative term associated to $\varphi^{\frac{ms+q-1}{ms}}$, instead of \eqref{eqdiff} we would write
\begin{equation*}
\varphi(t)'\leq c_1\varphi(t)+c_2\varphi^{\frac{3}{2}}(t)+c_3\varphi^3(t)\quad \textrm{for all}\quad t\in (0,t^*),
\end{equation*}
and the claim of the theorem would read
\[
t^*\geq \int_{\varphi(0)}^{+\infty}  \dfrac{d\tau}{c_1\tau+c_2\tau^{\frac{3}{2}}+c_5\tau^3},
\]
being in this case the last integral convergent but not explicitly computable.
\end{remark}
Through some straightforward manipulations we can prove that the previous theorem holds even in 2-dimensional settings, precisely as established in this 
\begin{theorem}\label{TheoremBlowUpDifferencePower2}  
Let $\Omega$ be a domain of $\mathbb{R}^2$ satisfying Assumption \ref{Assumption1}. Then, under the remaining assumptions of Theorem \ref{TheoremBlowUpDifferencePower}, 
there exist computable constants $\bar{c}_1>0$ and $\bar{c}_5>0$ such that for $\varphi(0)=\io u_0^{m(p-1)}>0$
 \begin{equation*} 
t^{*}    \geq \frac{1}{\bar{c}_1}\log \left(1+\frac{\bar{c}_1}{\bar{c}_5}\varphi^{-1}(0)\right)=:T.
 \end{equation*}
 In particular $I=[0,t^*)\supseteq [0,T).$
\begin{proof}
Similarly to what done throughout the proof of Theorem \ref{TheoremBlowUpDifferencePower}, let us rely on
\eqref{evolutionPhi4} and \eqref{Inequ_v^3_2_Precedente}. Conversely, in order to estimate $\io u^{\frac{3}{2}ms}$ we have now to refrain from using \eqref{Inequ_v^3} but \eqref{Inequ_v^3Dim2}, which yields for $\bar{\epsilon}>0$
\begin{equation*}
\begin{split}
\io u^{\frac{3}{2}ms}&\leq \dfrac{\sqrt{2}}{2\rho_0}\Big(\io u^{ms}\Big)^{\frac{3}{2}} +\dfrac{\sqrt{2}\left(d+\rho_0\right)}{4\rho_0\bar{\epsilon}_2^2}\Big(\io u^{ms}\Big)^2
\\ & \quad +\dfrac{\sqrt{2}\left(d+\rho_0\right)\bar{\epsilon}_2^2}{4\rho_0}\io \lvert\nabla u^{\frac{ms}{2}}\rvert^2\quad \textrm{for all}\quad t\in (0,t^*).
\end{split}
\end{equation*}
Subsequently, for $\varphi(t)=\varphi:=\io u^{ms}$ on $t\in (0,t^*)$, manipulations of the previous bound in conjunction with \eqref{evolutionPhi4} and \eqref{Inequ_v^3_2_Precedente} provide computable positive  constants $\bar{c}_1$, $\bar{c}_2=\bar{c}_2(\epsilon_1)$, $\bar{c}_3=\bar{c}_3(\bar{\epsilon}_2)$ (which we omit to calculate) and 
\begin{equation*}
\bar{c}_4=\bar{c}_4(\epsilon_1,\bar{\epsilon}_2)=\dfrac{\sqrt{2}}{4 \rho_0}(d +\rho_0){\bar{\epsilon}_2^2}\Big(2 a s|\Omega|+\dfrac{3m^2 s k}{2\rho_0}+\dfrac{5m^3s^2kd}{8 \epsilon_1\rho_0}\Big)
+\dfrac{5mdk\epsilon_1}{2\rho_0}-\dfrac{c}{ms}\,,
\end{equation*}
with the property that the following is ensured: 
\begin{equation*}
\varphi'\leq \bar{c}_1\varphi+\bar{c}_2\varphi^{\frac{3}{2}}+\bar{c}_3\varphi^2+\bar{c}_4\io \lvert\nabla u^{\frac{ms}{2}}\rvert^2 -bms\lvert\Omega\rvert^{\frac{1-q}{ms}}\varphi^{\frac{ms+q-1}{ms}}\quad \textrm{on}\quad (0,t^*).
\end{equation*}
As before, for any fixed $0<\epsilon_1<\frac{2\rho_0 c}{5m^2sdk}$ there exists $\bar{\epsilon}_2>0$ such that   $\bar{c}_4\leq 0$ and henceforth 
\begin{equation}\label{VariationPhiTwoD}
\varphi'\leq \bar{c}_1\varphi+\bar{c}_2\varphi^{\frac{3}{2}}+\bar{c}_3\varphi^2-bms\lvert\Omega\rvert^{\frac{1-q}{ms}}\varphi^{\frac{ms+q-1}{ms}}\quad \textrm{on}\quad (0,t^*).
\end{equation}
Using H\"{o}lder's and Young's inequality, we can estimate on $(0,t^*)$ the term involving $\varphi^{\frac{3}{2}}$  by means of a combination of $\varphi^{\frac{ms+q-1}{ms}}$ and $\varphi^2$, precisely obtaining 
\begin{equation*}
\varphi^{\frac{3}{2}}\leq \dfrac{ms}{2ms-2q+2}\bar{\epsilon}_3^{\frac{2ms-2q+2}{ms}}\varphi^{\frac{ms+q-1}{ms}}+\dfrac{ms-2q+2}{2ms-2q+2}\bar{\epsilon}_3^{-\frac{2ms-2q+2}{ms-2q+2}}\varphi^2.
\end{equation*}
In this way, expression \eqref{VariationPhiTwoD} reads
\begin{equation*}
\begin{split}
\varphi'\leq & \bar{c}_1\varphi+\left(\bar{c_2}\dfrac{ms}{2ms-2q+2}\bar{\epsilon}_3^{\frac{2ms-2q+2}{ms}}-bms\lvert\Omega\rvert^{\frac{1-q}{ms}}\right)\varphi^{\frac{ms+q-1}{ms}}\\&
+\Big(\bar{c}_3+\bar{c}_2\dfrac{ms-2q+2}{2ms-2q+2}\bar{\epsilon}_3^{-\frac{2ms-2q+2}{ms-2q+2}}\Big)\varphi^2\quad \textrm{for all}\quad t\in(0,t^*),
\end{split}
\end{equation*}
and for $\bar{\epsilon}_3=\left[\frac{2ms-2q+2}{\bar{c}_2}b\lvert\Omega\rvert^{\frac{1-q}{ms}}\right]^{\frac{ms}{2ms-2q+2}}$ and  $\bar{c}_5=\frac{ms-2q+2}{2ms-2q+2}\bar{c}_2\bar{\epsilon}_3^{-\frac{2ms-2q+2}{ms-2q+2}}+\bar{c}_3$ we finally have
\begin{equation}\label{eqdiff2}
\varphi(t)'\leq \bar{c}_1\varphi(t)+\bar{c}_5\varphi^2(t) \quad \textrm{for all}\quad t\in(0,t^*).
\end{equation}
As a consequence of these operations, and reasoning as in Theorem \ref{TheoremBlowUpDifferencePower}, our claim is given since
\begin{equation*}\label{Lower_Dirichlet}
\begin{split}
t^*\geq \int_{\varphi(0)}^{+\infty}  \dfrac{d\tau}{\bar{c}_1\tau+\bar{c}_5\tau^2} 
= \dfrac{1}{\bar{c}_1}\log \left(1+\dfrac{\bar{c}_1}{\bar{c}_5}\varphi^{-1}(0)\right).
\end{split}
\end{equation*}
\end{proof}
\end{theorem}
\begin{remark}
In line with Remark \ref{RemarkExplicitTStar}, but unlike its conclusion, if in \eqref{VariationPhiTwoD} we neglected the last negative part, instead of \eqref{eqdiff2} we would have
\begin{equation*}
\varphi(t)'\leq \bar{c}_1\varphi(t)+\bar{c}_2\varphi^{\frac{3}{2}}(t)+\bar{c}_3\varphi^2(t)\quad \textrm{on}\quad (0,t^*),
\end{equation*}
and the claim of the theorem would read
\[
t^*\geq \int_{\varphi(0)}^{+\infty}  \dfrac{d\tau}{\bar{c}_1\tau+\bar{c}_2\tau^{\frac{3}{2}}+\bar{c}_3\tau^2} =:T,
\]
being the last integral also explicitly computable if (and only if) $\Upsilon:=4\bar{c}_1\bar{c}_3-{\bar{c}_2}^2\geq 0$. More exactly, properties of inverse hyperbolic functions give  
\[
T=
\begin{cases}
-\frac{\bar{c}_2\pi-2\bar{c}_2 \arctan\left(\frac{\bar{c}_2+2 \bar{c}_3\sqrt{\varphi(0)}}{\sqrt{{-\bar{c}_2}^2+4 \bar{c}_1 \bar{c}_3}}\right)}{\bar{c}_1\sqrt{{-\bar{c}_2}^2+4 \bar{c}_1 \bar{c}_3}}-\frac{1}{\bar{c}_1}\log \big(\frac{\bar{c}_3\varphi(0)}{\bar{c}_1+\bar{c}_2 \sqrt{\varphi(0)}+\bar{c}_3 \varphi(0)}\big) & \textrm{if} \quad \Upsilon> 0,\\
\frac{-2}{\sqrt{c_1}(\sqrt{c_1}+\sqrt{c_1 \varphi(0)})} - \frac{1}{\bar{c}_1} \log\Big(\frac{\varphi(0)}{\big(\sqrt{\bar{c_1}}+\sqrt{\varphi(0)} \sqrt{\bar{c}_3}\big)^2}\Big)+ \frac{1}{\bar{c}_1}\log\left(\frac{1}{\bar{c}_3}\right) & \textrm{if} \quad \Upsilon= 0.
\end{cases}
\]
\end{remark}
  \subsection{A criterion for global existence}
In the last result, we are interested to examine the opposite situation described in Theorems \ref{TheoremBlowUpDifferencePower} and \ref{TheoremBlowUpDifferencePower2}. More exactly, we establish that when the effect of the source (coefficient $p$) is enough stronger than that of the diffusion (coefficient $m$) but  weaker than the one of the dampening (coefficient $q$), if a double stabilizing effect from the diffusion and the absorption somehow surpasses the same action of the source, system \eqref{General_Problem} does not suffer from blow-up phenomena, even for arbitrary large initial data $u_0({\bf x})$ and in presence of an incoming flow of the population $u$. 
\begin{theorem}\label{TheoremGlobalferencePower} 
Let $\Omega$ be a domain of $\R^N$, $N\geq 1,$ satisfying Assumption \ref{Assumption1}. Moreover, for $a, b, c, k > 0$, $q>p>m>1$, $2p<m+q$ let $0\leq g(\xi)\leq k\xi^{p-m+1}$, with $\xi \geq 0$.
If  $u$ is a classical solution, in the sense of Definition \ref{DefiSolution}, to \eqref{General_Problem} emanating from a positive initial data $u_0:=u_0({\bf x})\in C^{2+\alpha}(\Omega)\cap C^1(\bar{\Omega})$, for some $0<\alpha<1$ and such that  $\nabla u_0 \cdot\boldsymbol\nu=g(u_0)$ on $\partial \Omega$, then $t^*=\infty$, or equivalently $I=[0,\infty)$. 
\begin{proof}
If $u$ is a positive  classical solution of \eqref{General_Problem} defined in $\Omega \times (0,t^*)$ and satisfying $u_{ \boldsymbol\nu}=g(u)$ on $\partial \Omega$, by differentiating $\io u^{2}$ we derive 
\begin{equation}\label{diff}
\begin{split}
\frac{d}{dt}\io u^{2}  &=2\io  u \left(\Delta u^m+a\io u^p-bu^q-c\lvert\nabla\sqrt{u}\rvert^2\right) \\ & 
\leq 2mk\iob u^{p +1} -2m\io u^{m-1}\lvert \nabla u\rvert^2  +2a\lvert\Omega\rvert\io u^{p+1} -2b\io u^{q+1} \\&
-2c\io u\lvert\nabla\sqrt{u}\rvert^2\quad \textrm{for all}\quad t\in (0,t^*),
\end{split}
\end{equation}
where we have employed the following bound, consequence of the H\"{o}lder inequality:
\begin{equation*}
\io u\io u^p \leq \left(\io u^{p+1}\right)^{\frac{1}{p+1}}\lvert \Omega\rvert^{\frac{p}{p+1}}\left(\io u^{p+1}\right)^{\frac{p}{p+1}}\lvert\Omega\rvert^{\frac{1}{p+1}}=\lvert\Omega\rvert \io u^{p+1} \quad \textrm{on}\;(0,t^*)\,.
\end{equation*}
An application of \eqref{SobolevTypeInequBoundary} and the identity $
\dfrac{m+1}{2}u^p\lvert\nabla u\rvert =u^{\frac{2p-m+1}{2}}\lvert\nabla u^{\frac{m+1}{2}}\rvert$
provide
\begin{equation}\label{diff00}
\iob u^{p +1} \leq \dfrac{N}{\rho_0}\io u^{p +1} +\dfrac{2d\left(p+1\right)}{\rho_0\left(m+1\right)}\io u^{\frac{2p -m+1}{2}}\lvert\nabla u^{\frac{m+1}{2}}\rvert \quad\textrm{on}\;(0,t^*)\,.
\end{equation}
Now, by considering that
\begin{equation*}
-2m\io u^{m-1}\lvert\nabla u\rvert^2 =-\dfrac{8m}{\left(m+1\right)^2} \io \lvert \nabla u^{\frac{m+1}{2}}\rvert^2 \quad \textrm{on}\quad (0,t^*),
\end{equation*}
relation \eqref{diff} becomes by virtue of \eqref{diff00}
\begin{equation}\label{diff2}
\begin{split}
\frac{d}{dt}\io u^{2} \leq & \dfrac{2kNm}{\rho_0}\io u^{p+1}+\dfrac{4kdm\left(p+1\right)}{\rho_0\left(m+1\right)}\io u^{\frac{2p-m+1}{2}}\lvert\nabla u^{\frac{m+1}{2}}\rvert\\&
-\dfrac{8m}{\left(m+1\right)^2}\io \lvert\nabla u^{\frac{m+1}{2}}\rvert^2+2a\lvert\Omega\rvert\io u^{p+1} \\&-2b \io u^{q+1}\quad \textrm{for all}\;t\in (0,t^*)\,, 
\end{split}
\end{equation}
where, evidently, we have neglected the nonpositive term $-2c \io u\lvert\nabla\sqrt{u}\rvert^2$.

Additionally, from the Young inequality we obtain that for any $\sigma>0$
\begin{equation*}
\io u^{\frac{2p-m+1}{2}}\lvert\nabla u^{\frac{m+1}{2}}\rvert\leq \dfrac{\sigma}{2}\io u^{2p-m+1} +\dfrac{1}{2\sigma}\io \lvert\nabla u^{\frac{m+1}{2}}\rvert^2\quad \textrm{on}\;(0,t^*)\,,
\end{equation*}
so that fixing $\sigma=\frac{kd\left(p+1\right)\left(m+1\right)}{4\rho_0}$ this expression holds 
\begin{equation*}
\begin{split}
\dfrac{4kd\left(p+1\right)}{\rho_0\left(m+1\right)}m & \io u^{\frac{2p-m+1}{2}}\lvert\nabla u^{\frac{m+1}{2}}\rvert\leq  \dfrac{8m}{\left(m+1\right)^2}\sigma^2\io u^{2p-m+1}\\& \quad +\dfrac{8m}{\left(m+1\right)^2}\io \lvert\nabla u^{\frac{m+1}{2}}\rvert^2 \quad \textrm{on}\quad (0,t^*).
\end{split}
\end{equation*}
Combining this gained bound with \eqref{diff2}, we get
\begin{equation}\label{diff3}
\begin{split}
\frac{d}{dt}\io u^{2} &\leq \dfrac{2kNm}{\rho_0} \io u^{p+1}+\dfrac{8m}{\left(m+1\right)^2}\sigma^2\io u^{2p-m+1}\\&
+2a\lvert\Omega\rvert\io u^{p+1} -2b\io u^{q+1} \quad \textrm{on}\;(0,t^*)\,.
\end{split}
\end{equation}
Since $\frac{(q+1)(p-m)}{q-p}+\frac{(p+1)(q-2p+m)}{q-p}=2p-m+1$, Young's inequality produces for any $\epsilon>0$
\begin{equation}\label{young2}
\io u^{2p-m+1} \leq \left(1-\alpha\right)\epsilon\io u^{q+1}+\alpha\epsilon^{\frac{\alpha-1}{\alpha}}\io u^{p+1} \quad \textrm{on}\;(0,t^*)\,,
\end{equation}
where $0<\alpha =\frac{q+m-2p}{q-p}<1$ in view of the hypothesis $q>p>m$ and $2p<m+q$. Subsequently, by plugging \eqref{young2} into \eqref{diff3}, we obtain
\begin{equation}\label{EvolutionIntUsquare}
\begin{split}
\frac{d}{dt}\io u^2 &\leq \left(\dfrac{2kNm}{\rho_0}+2a\lvert\Omega\rvert+\dfrac{8m\sigma^2\alpha}{\left(m+1\right)^2}\epsilon^{\frac{\alpha-1}{\alpha}}\right) \io u^{p+1} \\ & 
\quad +\left(\dfrac{8m\sigma^2}{\left(m+1\right)^2}\left(1-\alpha\right)\epsilon-2b\right)\io u^{q+1}\\&
= M_1\io u^{p+1} -M_2\io u^{q+1} \quad \textrm{on}\;(0,t^*)\,,
\end{split}
\end{equation}
where 
\begin{equation*}
\begin{split}
M_1=\dfrac{2kNm}{\rho_0}+2a\lvert\Omega\rvert+\dfrac{8m\sigma^2\alpha}{(m+1)^2}\epsilon^{\frac{\alpha-1}{\alpha}}, \quad M_2=2b-\dfrac{8m\sigma^2\epsilon(1-\alpha)}{(m+1)^2}\,.
\end{split}
\end{equation*}
Now we let $\epsilon >0$ sufficiently small as to ensure $M_2>0$ (as an example we see that for $\epsilon=\frac{b(m+1)^2}{8m\sigma^2(1-\alpha)}$ we have $M_2=b$). Thereafter, since the Hölder inequality (recall $p <q$) gives 
\begin{equation*}
\io u^{p+1}\leq \left(\io u^{q+1}\right)^{\frac{p+1}{q+1}}\lvert\Omega\rvert^{\frac{q-p}{q+1}} \quad \textrm{on}\;(0,t^*),
\end{equation*}
we deduce from \eqref{EvolutionIntUsquare} that 
\begin{equation}\label{EvolutionUSquareAlmostFinal}
\begin{split}
\psi'\leq M_1\left(\io u^{q+1}\right)^{\frac{p+1}{q+1}}\left[\lvert\Omega\rvert^{\frac{q-p}{q+1}}-\frac{M_2}{M_1}\left(\io u^{q+1}\right)^{\frac{q-p}{q+1}}\right]\quad \textrm{on}\;(0,t^*),
\end{split}
\end{equation}
where we introduced $\psi(t)=\psi:=\io u^2$ for all $t\in(0,t^*)$. 
In order to establish an absorptive differential inequality for $\psi$, we use again Hölder's inequality to observe that
\begin{equation*}
\psi \leq \Big(\io u^{q+1}\Big)^\frac{2}{q+1}\lvert \Omega \rvert^\frac{q-1}{q+1}\Rightarrow -\Big(\io u^{q+1}\Big)^\frac{q-p}{q+1}\leq -\psi^\frac{q-p}{2}\lvert \Omega \rvert^\frac{(1-q)(q-p)}{2(q+1)} \quad \textrm{on}\quad (0,t^*),
\end{equation*}
so that in view of \eqref{EvolutionUSquareAlmostFinal} we arrive at this initial value problem,
\begin{equation*}
\begin{cases}
\psi'(t)\leq M_1\left(\io u^{q+1}\right)^{\frac{p+1}{q+1}}\left[\lvert\Omega\rvert^{\frac{q-p}{q+1}}-\frac{M_2}{M_1}\lvert \Omega \rvert^\frac{(1-q)(q-p)}{2(q+1)}\psi^\frac{q-p}{2}(t)\right]\quad \textrm{on}\;(0,t^*),\\
\psi(0)=\io u_0^2>0,
\end{cases}
\end{equation*}
from which is guaranteed that 
\begin{equation}\label{UniformBoundPsi}
\psi (t) \leq C:=\max\Bigg\{\io u_0^2,\lvert \Omega \rvert \bigg(\frac{M_1}{M_2}\bigg)^\frac{2}{q-p}\Bigg\} \quad \textrm{for all}\quad t\in(0,t^*).
\end{equation}
Finally, well know extension results for ODE's  with locally Lipschitz continuous right side (see, for instance, \cite{grant2014theoryODE}),  show that $t^*=\infty$; indeed, if $t^*$ were finite, $\psi (t) \nearrow +\infty$ as $t\searrow t^*$ and it would contradict \eqref{UniformBoundPsi}. 
\end{proof}
\end{theorem}
\subsubsection*{Acknowledgements}
The authors thank professor Stella Piro Vernier for her fruitful remarks and indications. MM and GV are members of the Gruppo Nazionale per l'Analisi Matematica, la Probabilit\`a e le loro Applicazioni (GNAMPA) of the Istituto Na\-zio\-na\-le di Alta Matematica (INdAM) and are partially supported by the research project \textit{Integro-differential Equations and Non-Local Problems}, funded by Fondazione di Sardegna (2017).

\end{document}